\documentclass[12pt,twoside]{amsart}
\usepackage{amssymb}
\usepackage{amscd}
\usepackage{longtable}

\title[Toric morphisms with anti-nef canonical divisors]
{Three-dimensional toric morphisms 
with anti-nef canonical divisors} 
\author{Hiroshi Sato} 
\thanks{The author is partly supported by the 
Grant-in-Aid for JSPS Fellows, The Ministry of 
Education, Science, Sports and Culture, Japan.}
\address{
Osaka City University Advanced Mathematical Institute, 
3-3-138 Sugimoto, Sumiyoshi-ku, Osaka 558-8585, 
Japan}
\email{hirosato@sci.osaka-cu.ac.jp}

\newcommand{\Hom}[0]{{\operatorname{Hom}}}

\newcommand{\Conv}[0]{{\operatorname{Conv}}}
\newcommand{\Cone}[0]{{\operatorname{Cone}}}
\newcommand{\G}[0]{{\operatorname{G}}}
\newcommand{\Int}[0]{{\operatorname{Int}}}
\newtheorem{thm}{Theorem}[section]
\newtheorem{lem}[thm]{Lemma}

\newtheorem{prop}[thm]{Proposition}
\newtheorem{lem-def}[thm]{Definition-Lemma}

\theoremstyle{definition}
\newtheorem{ex}[thm]{Example}
\newtheorem{defn}[thm]{Definition}
\newtheorem{prob}[thm]{Problem}
\newtheorem{rem}[thm]{Remark}
\newtheorem*{ack}{Acknowledgments}       
 
\begin{document}
\bibliographystyle{amsalpha+}

\begin{abstract}
In this paper, we classify projective toric birational 
morphisms from Gorenstein toric $3$-folds onto 
the $3$-dimensional affine space with 
relatively ample anti-canonical divisors. 
\end{abstract}

\maketitle
\tableofcontents

\section{Introduction}\label{intro}

\thispagestyle{empty}

The classification of toric (weak) Fano varieties is an 
important problem in the toric geometry. 
It is well-known that the number of Gorenstein toric Fano 
$d$-folds up to isomorphisms is {\em finite} for any $d$ and 
Gorenstein toric Fano $d$-folds are classified for $d\leq 4$ 
(see \cite{batyrev3}, \cite{watanabewatanabe}, 
\cite{batyrev4} and \cite{sato} for smooth cases, 
and see \cite{koelman1}, 
\cite{ks3} and \cite{ks4} for Gorenstein cases). 
In this paper, as a local version of 
this problem, we consider the following problem which was 
posed by Watanabe:

\begin{prob}\label{afureru}
Classify all the projective toric birational 
morphisms $f:X\to{\mathbb A}^3$ such that $X$ is smooth and 
$K_X$ is a relatively anti-nef divisor. 
\end{prob}
By considering this problem up to {\em flops}, Problem \ref{afureru} 
is essentially equivalent to the following problem:

\begin{prob}\label{oogesane}
Classify all the projective toric birational 
morphisms $f:X\to{\mathbb A}^3$ such that $X$ is Gorenstein and 
$-K_X$ is a relatively ample divisor. 
\end{prob}

We solve Problem \ref{oogesane} under one more assumption that 
the image of the exceptional set of $f$ is the origin 
(see Definition \ref{iinokana}). 
One of the difficulties of this problem is that the 
number of such morphisms up to isomorphisms is {\em infinite} 
(see Remark \ref{akitaken}), while 
the number of genuine 
Gorenstein toric Fano $3$-folds up to isomorphisms 
is finite as above. However, we obtain the explicit 
description of the classification. There are three classes of 
such morphisms (see Theorem \ref{babyg}). 

The content of this paper is as follows: In Section \ref{miyagiken}, 
we give the definition of {\em local toric Fano} varieties. Then, 
we translate the properties of local toric Fano varieties into 
the properties of fans. Moreover, we consider the $2$-dimensional 
version of Problem \ref{oogesane}. Section \ref{daiei} is 
devoted to the classification of Gorenstein local toric Fano 
$3$-folds, which is a partial answer to Problem \ref{oogesane}. 
The explicit classified list of Gorenstein local toric Fano 
$3$-folds will be given.

\begin{ack}
The author would like to thank Professor Kei-ichi Watanabe for 
introducing him to this problem 
and giving useful comments. 
He also would like to thank Professors 
Osamu Fujino and Masataka Tomari 
for advice and encouragement.
\end{ack}

\section{The definition of local toric Fano varieties}
\label{miyagiken}

First of all, we prepare the notation. For the fundamental 
properties of the toric geometry, see \cite{fulton} and \cite{oda}. 

Let $N:={\mathbb Z}^d$ and $M:=\Hom_{\mathbb Z}(N,{\mathbb Z})$ 
the dual group. 
The natural pairing $\langle\ ,\ \rangle :M\times N \to
{\mathbb Z}$ 
is extended to a bilinear form 
$\langle\ ,\ \rangle :M_{\mathbb R}\times N_{\mathbb R}\to{\mathbb R}$, 
where $M_{\mathbb R}:=M\otimes_{\mathbb Z}{\mathbb R},\ N_{\mathbb R}
:=N\otimes_{\mathbb Z}{\mathbb R}$. Let $X=X_\Sigma$ be 
a Gorenstein toric $d$-fold associated to a fan $\Sigma$ in $N$. 
Let 
$\G(\Sigma)$ be the set of primitive generators 
of $1$-dimensional cones in $\Sigma$, and put 
$\G(\sigma):=\sigma\cap\G(\Sigma)$ for a cone 
$\sigma\in\Sigma$. 
For $\sigma\in\Sigma$, we put 
$F_{\sigma}:=\Conv(\G(\sigma))$ and 
$\Gamma_{\sigma}:=\Conv(\G(\sigma)\cup\{0\})$, where 
$\Conv(T)$ is the convex hull for any subset $T\subset N_{\mathbb R}$. 
For a $d$-dimensional cone $\sigma\in\Sigma$, there exists the 
unique element $m_{\sigma}\in M$ such that the hyperplane 
$H_{\sigma}\subset N_{\mathbb R}$ 
generated by $\G(\sigma)$ is defined by $m_\sigma$, that is,
$$H_\sigma=\left\{n\in N_{\mathbb R}\,\left|\,
\langle m_{\sigma},n\rangle=-1\right.\right\}.$$
We will use these notation throughout this paper. 

\begin{defn}\label{iinokana}
Let $X=X_{\Sigma}$ be a Gorenstein toric $d$-fold. 
$X$ is called a local toric 
{\em Fano} (resp. {\em weak Fano}) $d$-fold if 
there exists a projective birational toric morphism 
$f:X\to{\mathbb A}^d$, where ${\mathbb A}^d$ is 
the $d$-dimensional affine space, 
$-K_X$ is $f$-ample (resp. $f$-nef) and the image of 
the exceptional set of $f$ is the origin. 
\end{defn}

\begin{rem}
Smooth local toric Fano $3$-folds and $4$-folds 
were studied in \cite{sato} and \cite{casagrande}, respectively. 
In these papers, the final assumption that the image of 
the exceptional set of $f$ is the origin was {\em not} necessary.
\end{rem}

Let $\{e_1,\ldots,e_d\}$ be the standard basis for $N$, 
$\sigma_0:=\Cone(e_1,\ldots,e_d)$ and 
$\Sigma_0$ the fan whose cones are all the faces of $\sigma_0$. 
Then, the fan corresponding to ${\mathbb A}^d$ is $\Sigma_0$. 
The following is a characterization of local toric Fano $d$-folds 
using the properties of fans. 
The proof is similar as the case of genuine toric Fano $d$-folds. 

\begin{prop}\label{charfano}
Let $X=X_{\Sigma}$ be a Gorenstein toric $d$-fold. $X$ is a local 
toric weak Fano $d$-fold if and only if $|\Sigma|=\sigma_0$, 
$\tau\in\Sigma$ for any $\tau\in\Sigma_0\setminus\{\sigma_0\}$ and 
$$\Gamma_{\Sigma}:=\bigcup_{\sigma\in\Sigma}\Gamma_{\sigma}$$ is 
a convex subset in $N_{\mathbb R}$. 
On the other hand, 
$X$ is a local toric 
Fano $d$-fold if and only if $|\Sigma|=\sigma_0$, 
$\tau\in\Sigma$ for any $\tau\in\Sigma_0\setminus\{\sigma_0\}$ 
and 
$\Gamma_{\Sigma}$ is 
a strictly convex subset in $N_{\mathbb R}$, that is, 
$\Gamma_{\Sigma}$ is convex and 
for any $d$-dimensional cone $\sigma\in\Sigma$ and for any 
$x\in\Gamma_{\Sigma}$, $\langle m_{\sigma},x\rangle=-1$ 
implies $x\in F_{\sigma}$.
\end{prop}

The main purpose of this paper is to give the classification 
of smooth local toric weak Fano $3$-folds up to 
isomorphisms, which was posed by Watanabe. 
However, this problem seems very complicated because  
there may exist many flopping contractions. 
So, we consider this problem up to flops, that is, 
we do not classify smooth local toric weak Fano $3$-folds 
but Gorenstein local toric Fano $3$-folds. 
Then, in principle, we can obtain the classification 
of smooth local toric weak Fano $3$-folds 
by toric crepant resolutions (see \cite{odapark}). 

As the first step for our main result, 
we consider the classification of Gorenstein local toric 
Fano {\em surfaces}. In this case, there exists a 
one-to-one corresponding between Gorenstein local toric 
Fano surfaces and smooth local toric weak 
Fano surfaces. Namely, we obtain 
a Gorenstein local toric Fano surface 
by the anti-canonical morhism of 
a smooth local toric weak Fano surface, while 
we obtain a smooth local toric weak Fano surface 
by the toric crepant resolution of 
a Gorenstein local toric Fano surface. 

\begin{ex}[The classification of smooth 
local toric weak Fano surfaces]\label{tontoro}
Let $X=X_{\Sigma}$ be a smooth local toric weak 
Fano surfaces. 
The morphism $f:X\to{\mathbb A}^2$ is factored into 
a sequence of toric blow-ups:
$$X=:X_n\to X_{n-1}\to \cdots\to X_1\to X_0:={\mathbb A}^2.$$ 
If $n\neq 0$, then ${}^t(1,1)\in\G(\Sigma)$ 
because $X_1\to X_0$ is a toric blow-up. Similarly, 
${}^t(2,1)$ or ${}^t(1,2)$ is contained in $\G(\Sigma)$ 
if $n\geq 2$. 
Now, suppose that ${}^t(2,1)\in\G(\Sigma)$. If $\G(\Sigma)$ 
contains an element ${}^t(s,t)\in N$ such that 
$s\geq1$ and $t\geq 2$, then 
$$\Conv\left(0,
\begin{pmatrix}
1\cr 0
\end{pmatrix},
\begin{pmatrix}
0\cr 1
\end{pmatrix},
\begin{pmatrix}
2\cr 1
\end{pmatrix},
\begin{pmatrix}
s\cr t
\end{pmatrix}
\right)$$
contains ${}^t(1,1)$ in its interior. 
This contradicts the Gorenstein condition. 
Thus, one can easily check that $X$ is uniquely 
determined by $n$, that is, 
$$\G(\Sigma)=\left\{
\begin{pmatrix}
1\cr 0
\end{pmatrix},
\begin{pmatrix}
0\cr 1
\end{pmatrix},
\begin{pmatrix}
1\cr 1
\end{pmatrix},\ldots,
\begin{pmatrix}
n\cr 1
\end{pmatrix}
\right\}.$$
In Figure \ref{fanosurface}, we describe $\Gamma_{\Sigma}$. 
For the corresponding Gorenstein local toric 
Fano surface, we obtain the same convex set. 

\bigskip

\begin{figure}[hbtp]
\centering
\begin{picture}(230,120)
\put(50,0){\vector(0,1){120}}
\put(0,50){\vector(1,0){230}}
\put(65,20){
$\begin{pmatrix}
1\cr 0
\end{pmatrix}$
}
\put(15,80){
$\begin{pmatrix}
0\cr 1
\end{pmatrix}$
}
\put(40,38){$0$}
\put(70,62){$\Gamma_{\Sigma}$}

\put(80,50){\circle*{5}}
\put(50,50){\circle*{5}}
\put(50,80){\circle*{5}}
\thicklines
\put(50,80){\line(1,0){150}}
\put(80,50){\line(4,1){120}}
\put(50,50){\line(0,1){30}}
\put(50,50){\line(1,0){30}}
\thinlines

\put(200,80){\circle*{5}}
\put(200,95){
$\begin{pmatrix}
n\cr 1
\end{pmatrix}$
}
\end{picture}
\caption{}
\label{fanosurface}
\end{figure}

\end{ex}

\begin{rem}\label{akitaken}
As Example \ref{tontoro}, there exist {\em infinitely} many 
Gorenstein local toric Fano $d$-folds for $d\geq 2$.
\end{rem}

\section{Gorenstein local toric Fano 3-folds}
\label{daiei}

In this section, we classify Gorenstein local toric Fano $3$-folds. 

First of all, we give examples of 
Gorenstein local toric Fano $3$-folds. 
This large class of Gorenstein local toric Fano $3$-folds 
is important for our classification. 

\begin{ex}\label{pkarasada}
Let $(Z_1,Z_2,Z_3)$ be the coordinate for $N_{\mathbb R}\cong{\mathbb R}^3$. 
Let $P$ be a convex polytope in the hyperplane 
$(Z_3=1)\subset N_{\mathbb R}$ which satisfies 
the following conditions:
\begin{enumerate}
\item $e_3$ is a vertex of $P$.
\item For any $x={}^t(a,b,1)\in P$, we have $a>0$ and $b>0$ unless 
$x=e_3$.
\item Put the vertices of $P$ be 
$$
x_0:=e_3=
\begin{pmatrix}
0\cr 0\cr 1
\end{pmatrix},
x_1=
\begin{pmatrix}
a_1\cr b_1\cr 1
\end{pmatrix},\ldots,
x_n=
\begin{pmatrix}
a_n\cr b_n\cr 1
\end{pmatrix}
$$
clockwise (see Figure \ref{polykara}). 
We put $x_{n+1}=x_0$ for convenience. 
Then, for any $1\leq i\leq n+1$, we have the following:
$$
\left\{
\begin{array}{rcl}
\frac{b_i-b_{i-1}}{a_i-a_{i-1}}\in{\mathbb Z} 
& \mbox{if} & a_i-a_{i-1}>0. \\
\frac{a_i-a_{i-1}}{b_i-b_{i-1}}\in{\mathbb Z} 
& \mbox{if} & a_i-a_{i-1}\leq 0. \\
\end{array}
\right.
$$
\end{enumerate}
For such a polytope $P$, we can define 
the fan $\Sigma_P$ whose maximal cones are 
\begin{enumerate}
\item $\Cone(P)$, 
\item $\Cone(e_2,x_{i-1},x_i)$ for $b_i-b_{i-1}\geq-(a_i-a_{i-1})$,
\item $\Cone(e_1,e_2,x_{i-1},x_i)$ for 
$b_i-b_{i-1}=-(a_i-a_{i-1})$ and 
\item $\Cone(e_1,x_{i-1},x_i)$ for $b_i-b_{i-1}\leq-(a_i-a_{i-1})$. 
\end{enumerate}
Then, the corresponding toric $3$-fold $X_P:=X_{\Sigma_P}$ is 
a Gorenstein local toric Fano $3$-fold. 
To confirm this, we may suffice to say that 
the maximal cones are Gorenstein. 
Indeed, the equations of hyperplanes 
generated by the primitive generators of the above maximal cones are 
$$(1)\ Z_3=1,$$
$$(2)\ \frac{b_i-b_{i-1}}{a_i-a_{i-1}} Z_1-Z_2+
\left(-1+b_i-\frac{b_i-b_{i-1}}{a_i-a_{i-1}} a_i\right)Z_3=-1,
$$
$$
(3)\ -Z_1-Z_2+(-1+a_i+b_i)Z_3=-1\ \mbox{and}
$$
$$(4)\ -Z_1+\frac{a_i-a_{i-1}}{b_i-b_{i-1}} Z_2+
\left(-1+a_i-\frac{a_i-a_{i-1}}{b_i-b_{i-1}} b_i\right)Z_3=-1,
$$
respectively. We remark that all the coefficients are integers. 

\bigskip

\begin{figure}[hbtp]
\centering
\begin{picture}(260,240)
\put(50,0){\vector(0,1){240}}
\put(0,50){\vector(1,0){260}}

\put(130,130){$P$}

\thicklines
\put(50,50){\circle*{5}}
\put(50,50){\line(1,3){30}}
\put(35,38){$x_0$}

\put(80,140){\circle*{5}}
\put(80,140){\line(1,2){30}}
\put(65,150){$x_1$}

\put(110,200){\circle*{5}}
\put(110,200){\line(1,0){60}}
\put(100,210){$x_2$}

\put(170,200){\circle*{5}}
\put(170,200){\line(1,-1){30}}
\put(200,170){\circle*{5}}
\put(200,170){\line(0,-1){60}}
\put(200,110){\circle*{5}}
\put(200,110){\line(-1,-1){30}}
\put(205,100){$x_{n-1}$}

\put(170,80){\circle*{5}}
\put(170,80){\line(-4,-1){120}}
\put(175,70){$x_n$}

\thinlines

\end{picture}
\caption{}
\label{polykara}
\end{figure}

\end{ex}

\begin{prop}\label{sakura}
Let $X=X_{\Sigma}$ be a 
Gorenstein local toric Fano $3$-fold. 
$\Sigma$ is associated to a polytope 
$P\subset (Z_3=1)\subset N_{\mathbb R}$ as in 
{\rm Example} $\ref{pkarasada}$ 
if and only if there exists a $3$-dimensional cone $\sigma\in\Sigma$ 
such that $H_{\sigma}=(Z_3=1)$.
\end{prop}

\begin{proof}
Suppose that there exists a $3$-dimensional cone $\sigma\in\Sigma$ 
such that $H_{\sigma}=(Z_3=1)$. 
For any $x\in\G(\Sigma)\setminus\G(\sigma)$, 
$x$ is contained in $(Z_3=0)$. So, $x$ is either $e_1$ or $e_2$ 
by Propsition \ref{charfano}. Thus, because of the Gorenstein 
condition, $X$ is associated to a polytope 
$P\subset (Z_3=1)\subset N_{\mathbb R}$ 
as in Example $\ref{pkarasada}$. Another implication is trivial. 
\end{proof}

Now, we start the classification. 
The starting point is the following simple lemma.

\begin{lem}\label{simplekey}
Let $f:X\to{\mathbb A}^d$ 
be a projective toric birational morphism from 
a smooth toric $d$-fold $X=X_{\Sigma}$. 
If the image of the exceptional set of $f$ is the origin, 
then either $f$ is an isomorphism or 
$e_1+\cdots+e_d\in\G(\Sigma).$

\end{lem}

\begin{proof}
$e_1+\cdots+e_d$ is contained in the relative interior of 
some cone $\sigma\in\Sigma$. If $\dim\sigma=d$, then $f$ is 
an isomorphism. So, suppose that $1<\dim\sigma<d$. 
Then, at least one element $x$ of $\G(\sigma)$ is 
not contained in $\{e_1,\ldots,e_d\}$. 
Thus, some $(d-1)$-dimensional cone in $\Sigma_0$ is 
not contained in $\Sigma$. This contradicts the fact that 
the image of the exceptional set of $f$ is the origin. 
Therefore, $\dim\sigma=1$.
\end{proof}

\renewcommand{\theenumi}{\Roman{enumi}}

Let $X=X_{\Sigma}$ be a Gorenstein local toric Fano $3$-fold and 
$\widetilde{X}\to X$ a toric crepant resolution with 
$\widetilde{X}=X_{\widetilde{\Sigma}}$. $\widetilde{X}$ is a 
smooth local toric weak Fano $3$-fold. Assume that $f$ is 
not an isomorphism. Then, by Lemma \ref{simplekey}, 
$x_0:=e_1+e_2+e_3\in\G(\widetilde{\Sigma})$. 
Thus, we have the four situations:
\begin{enumerate}
\item $\Cone(x_0,e_1)\in\widetilde{\Sigma}$, 
$\Cone(x_0,e_2)\in\widetilde{\Sigma}$ and 
$\Cone(x_0,e_3)\in\widetilde{\Sigma}$.
\item $\Cone(x_0,e_1)\not\in\widetilde{\Sigma}$, 
$\Cone(x_0,e_2)\in\widetilde{\Sigma}$ and 
$\Cone(x_0,e_3)\in\widetilde{\Sigma}$.
\item $\Cone(x_0,e_1)\not\in\widetilde{\Sigma}$, 
$\Cone(x_0,e_2)\not\in\widetilde{\Sigma}$ and 
$\Cone(x_0,e_3)\in\widetilde{\Sigma}$.
\item $\Cone(x_0,e_1)\not\in\widetilde{\Sigma}$, 
$\Cone(x_0,e_2)\not\in\widetilde{\Sigma}$ and 
$\Cone(x_0,e_3)\not\in\widetilde{\Sigma}$.
\end{enumerate}

\begin{lem}
If $\Cone(x_0,e_i)\not\in\widetilde{\Sigma}$, then 
$x_i:=x_0+e_i\in\G(\widetilde{\Sigma})$.
\end{lem}

\begin{proof}
We prove the lemma for the case $i=1$. 

Suppose that $x_1\not\in\G(\widetilde{\Sigma})$. 
Then, we have the equality 
$$x_1=
\begin{pmatrix}
2\cr 1\cr 1
\end{pmatrix}
=y_1+y_2$$
for $\{y_1,y_2\}\subset\G(\widetilde{\Sigma})$ 
which differs from $\{x_0,e_1\}$. 
Therefore, at least one of the $2$-dimensional cones in 
$\Sigma_0$ is not contained in $\widetilde{\Sigma}$. 
This is a contradiction.
\end{proof}

In Figure \ref{baaiwake}, 
we give the pictures of the above four situations 
by giving the sections of cones as usual. 

\begin{figure}[hbtp]
\centering
\begin{picture}(350,310)


\put(10,290){{\rm (I)}}

\put(30,180){\line(3,5){60}} 
\put(30,180){\line(1,0){120}}
\put(150,180){\line(-3,5){60}}

\put(30,180){\line(5,3){60}}
\put(150,180){\line(-5,3){60}}
\put(90,280){\line(0,-1){64}}

\put(90,216){\circle*{5}}
\put(30,180){\circle*{5}}
\put(150,180){\circle*{5}}
\put(90,280){\circle*{5}}

\put(85,201){$x_0$}
\put(15,170){$e_2$}
\put(155,170){$e_3$}
\put(85,290){$e_1$}


\put(180,290){{\rm (II)}}

\put(200,180){\line(3,5){60}}
\put(200,180){\line(1,0){120}}
\put(320,180){\line(-3,5){60}}

\put(200,180){\line(5,3){60}}
\put(320,180){\line(-5,3){60}}

\put(260,216){\circle*{5}}
\put(200,180){\circle*{5}}
\put(320,180){\circle*{5}}
\put(260,280){\circle*{5}}

\put(260,248){\circle*{5}}
\put(255,233){$x_1$}


\put(10,140){{\rm (III)}}

\put(30,30){\line(3,5){60}}
\put(30,30){\line(1,0){120}}
\put(150,30){\line(-3,5){60}}

\put(150,30){\line(-5,3){60}}

\put(90,66){\circle*{5}}
\put(30,30){\circle*{5}}
\put(150,30){\circle*{5}}
\put(90,130){\circle*{5}}

\put(90,98){\circle*{5}}
\put(60,48){\circle*{5}}
\put(70,43){$x_2$}


\put(180,140){{\rm (IV)}}

\put(200,30){\line(3,5){60}}
\put(200,30){\line(1,0){120}}
\put(320,30){\line(-3,5){60}}

\put(260,66){\circle*{5}}
\put(200,30){\circle*{5}}
\put(320,30){\circle*{5}}
\put(260,130){\circle*{5}}

\put(260,98){\circle*{5}}
\put(230,48){\circle*{5}}
\put(290,48){\circle*{5}}

\put(270,43){$x_3$}

\thicklines

\end{picture}
\caption{}
\label{baaiwake}
\end{figure}

We consider the classification 
of Gorenstein local toric Fano $3$-folds 
for these four cases.

\bigskip

(I) In this case, $f$ is factored into 
$$\widetilde{X}\to X_1\to {\mathbb A}^3,$$
where $X_1\to {\mathbb A}^3$ is the toric blow-up at the origin. 
Let $\Sigma_1$ be the fan corresponding to $X_1$. Then, 
the following holds.

\begin{lem}\label{tonnels}
At least two $3$-dimensional cones in $\Sigma_1$ 
are in $\widetilde{\Sigma}$.
\end{lem}

\begin{proof}
Suppose that $\Cone(x_0,e_1,e_2)$ and $\Cone(x_0,e_2,e_3)$ 
are not in $\widetilde{\Sigma}$. Then, by Lemma \ref{simplekey}, 
$y_1:=x_0+e_1+e_2$ and $y_2:=x_0+e_2+e_3$ 
are in $\G(\widetilde{\Sigma})$. 
However, we have $y_1+y_2=2x_0+e_1+2e_2+e_3=3x_0+e_2$, and 
this contradicts the convexity of $\Gamma_{\widetilde{\Sigma}}$. 
\end{proof}

We may assume that $\Cone(x_0,e_2,e_3)$ and $\Cone(x_0,e_3,e_1)$ 
are in $\widetilde{\Sigma}$ by Lemma \ref{tonnels}. The following 
proposition completes the classification for the case (I).

\begin{prop}\label{case1matome}
$$\Gamma_{\Sigma}=\Conv\left(0,e_1,e_2,e_3,y_m:=
\begin{pmatrix}
m\cr m\cr 1
\end{pmatrix}
\right)
$$
for some integer $m\geq 1$ $($see {\rm Figure $\ref{case1}$}$)$.
\end{prop}

\begin{proof}
If $\widetilde{X}\to X_1$ is not an isomorphism, 
then there exists $y\in\G(\widetilde{\Sigma})$ which is contained in 
the interior of $\Cone(e_1,e_2,x_0)$. So, let 
$$y:=ae_1+be_2+cx_0=
\begin{pmatrix}
a+c\cr b+c\cr c
\end{pmatrix}$$
for positive integers $a,b,c$. Suppose that $a>b$. 
Then, we have 
$$y+be_3=
\begin{pmatrix}
a+c\cr b+c\cr b+c
\end{pmatrix}
=(a-b)e_1+(b+c)x_0.$$
The inequality $1+b<c+a=(a-b)+(b+c)$ contradicts the convexity 
of $\Gamma_{\widetilde{\Sigma}}$. Therefore, $a=b$. Similarly, 
the equality
$$y+ae_3=
\begin{pmatrix}
a+c\cr a+c\cr a+c
\end{pmatrix}
=(a+c)x_0$$
implies that $c=1$. Finally, by putting 
$$m:=\max\left\{a\in{\mathbb Z}_{>0}\,\left|\,
\begin{pmatrix}
a\cr a\cr 1
\end{pmatrix}
\in\G(\widetilde{\Sigma})\right.\right\},
$$
we complete the proof.
\end{proof}

\bigskip

\begin{figure}[hbtp]
\centering
\begin{picture}(180,160)

\put(30,30){\line(3,5){60}}
\put(30,30){\line(1,0){120}}
\put(150,30){\line(-3,5){60}}

\put(30,30){\line(5,6){40}}
\put(150,30){\line(-5,3){80}}
\put(90,130){\line(-2,-5){20}}

\put(70,78){\circle*{5}}

\put(67,63){$y_m$}


\put(30,30){\circle*{5}}
\put(150,30){\circle*{5}}
\put(90,130){\circle*{5}}


\end{picture}
\caption{}
\label{case1}
\end{figure}

\bigskip

(II) This case is the main part of the classification. 

\begin{lem}
$\Cone(e_2,e_3,x_0)\in\widetilde{\Sigma}$.
\end{lem}

\begin{proof}
If $\Cone(e_2,e_3,x_0)\not\in\widetilde{\Sigma}$, 
then $e_2+e_3+x_0={}^t(1,2,2)\in\G(\widetilde{\Sigma})$ by 
Lemma \ref{simplekey}. 
$$x_1+
\begin{pmatrix}
1\cr 2\cr 2
\end{pmatrix}
=
\begin{pmatrix}
3\cr 3\cr 3
\end{pmatrix}
=3x_0$$
implies that $\Gamma_{\Sigma}$ is not convex. 
Thus, $\Cone(e_2,e_3,x_0)\in\widetilde{\Sigma}$.
\end{proof}

\renewcommand{\theenumi}{II\alph{enumi}}

By applying the classification of smooth 
local toric weak Fano surfaces (see Example \ref{tontoro}) 
for the $2$-dimensional smooth cone $\Cone(e_1,x_0)$, 
we have the following:

\begin{lem}
One of the following holds$:$
\begin{enumerate}
\item There exists a integer $l\geq 0$ such that
$$\G(\widetilde{\Sigma})\cap\Cone(e_1,x_0)
=\{x_0,x_1,x_1+e_1,\ldots,x_1+le_1,e_1\}.$$
\item $$\G(\widetilde{\Sigma})\cap\Cone(e_1,x_0)
=\{x_0,x_0+x_1,x_1,e_1\}.$$
\end{enumerate}
\end{lem}

\begin{proof}
For case (IIb), there does not exist element in 
$\G(\widetilde{\Sigma})\cap\Cone(e_1,x_0)$ other than 
$\{x_0,x_0+x_1,x_1,e_1\}$ 
because the hyperplane spaned by $e_2$, $e_3$ and $x_0+x_1$ 
contains $x_0$. 
\end{proof}

(IIa) Suppose that there is an element 
$z\in\G(\widetilde{\Sigma})\cap\Int\left(\Cone(e_1,e_2,x_0)\right)$, 
where $\Int(S)$ stands for the interior of $S$. 
Put $z:=ae_1+be_2+cx_0$ for positive integers 
$a,b,c\in{\mathbb Z}_{>0}$. 

First, suppose that $c=1$. Then, it is obvious that 
$\Sigma$ is associated to a polytope $P\subset(Z_3=1)$ 
which contains $e_3$, $x_0$, $x_1+le_1$ and $z$ 
as 
in Example \ref{pkarasada}.

So, suppose that $c\geq 2$. But, this is impossible because 
$\Conv(x_0,x_1)$ is contained 
in the boundary of $\Gamma_{\Sigma}$. Therefore, there 
exists no such $z$. 

We can apply these arguments for $\left(\Cone(e_1,e_3,x_0)\right)$ 
similarly. Thus, in this case, 
$\Sigma$ is associated to a polytope $P\subset(Z_3=1)$ as 
in Example \ref{pkarasada}. 

\bigskip

(IIb) First of all, we remark that $\Conv(e_2,e_3,x_0+x_1)$ and 
$\Conv(e_1,x_0+x_1)$ are contained in the boundary of 
$\Gamma_{\Sigma}$. The equations of hyperplanes spaned by 
$\{e_2,e_3,x_0+x_1\}$ and $\{e_1,e_3,x_0+x_1\}$ are 
$$Z_1-2Z_2+Z_3=1\mbox{ and }-Z_1+Z_2+Z_3=1,$$
respectively. 
Let $\tau$ be the unique $3$-dimensional cone in $\Sigma$ 
such that $e_1,e_2\in\tau\cap\G(\Sigma)$ and the equation 
of $H_{\tau}$
$$Z_1+Z_2+\alpha Z_3=1,$$
where $\alpha\in{\mathbb Z}$. 
Then, if there is an element 
$z\in\G(\widetilde{\Sigma})\cap
\Int\left(\Cone(e_1,e_2,x_0+x_1)\right)$, then 
$z$ is contained in the area
$$\{Z_1-2Z_2+Z_3\leq 1\}\cap\{-Z_1+Z_2+Z_3\leq 1\}\cap
\{Z_1+Z_2+\alpha Z_3\leq 1\}.$$
Any point $x$ on the line 
$$L=(Z_1-2Z_2+Z_3=1)\cap(-Z_1+Z_2+Z_3=1)$$
is described as 
$$x=
\begin{pmatrix}
3t-3\cr 2t-2\cr t
\end{pmatrix}
$$
for $t\in{\mathbb R}$. If $L$ and $H_{\tau}$ intersect 
at $t_0\in{\mathbb R}$, then we have 
$$\alpha=\frac{6}{t_0}-5.$$
Obviously, $t_{0}\geq 1$. Thus, we have $1\leq t_{0}\leq 6$, 
because $\alpha\in{\mathbb Z}$. Therefore, 
$z$ is contained in the area
$$\{Z_1-2Z_2+Z_3\leq 1\}\cap\{-Z_1+Z_2+Z_3\leq 1\}\cap
\{Z_1+Z_2-4Z_3\leq 1\}.$$

By considering the area 
$\Int\left(\Cone(e_1,e_3,x_0+x_1)\right)$ similarly, 
we have the following:

\begin{prop}
Put 
$$
z_1:=\begin{pmatrix}
15\cr 10\cr 6
\end{pmatrix}
\mbox{ and }
z_2:=
\begin{pmatrix}
15\cr 6\cr 10
\end{pmatrix}
.$$
Then, 
$\Gamma_{\Sigma}$ is contained in 
$$\Conv\left(
0,e_1,e_2,e_3,z_1
\right)\cup
\Conv\left(
0,e_1,e_2,e_3,z_2
\right).$$
This area is {\em not} 
a convex set 
$($see {\rm Figure \ref{mainnandabesa}}$)$. 
In particular, the number of 
Gorenstein local toric Fano $3$-folds 
is {\em finite} in this case.
\end{prop}

\begin{figure}[hbtp]
\centering
\begin{picture}(180,160)

\put(30,30){\line(3,5){60}}
\put(30,30){\line(1,0){120}}
\put(150,30){\line(-3,5){60}}

\put(150,30){\line(-5,6){40}}

\put(30,30){\line(5,6){40}}

\put(90,130){\line(-2,-5){20}}
\put(90,130){\line(2,-5){20}}
\put(90,130){\line(0,-1){64}}

\put(90,66){\circle*{5}}

\put(71,46){$x_0+x_1$}

\put(70,78){\circle*{5}}

\put(65,65){$z_1$}

\put(70,78){\line(5,-3){20}}
\put(110,78){\line(-5,-3){20}}

\put(110,78){\circle*{5}}

\put(107,65){$z_2$}



\put(30,30){\circle*{5}}
\put(150,30){\circle*{5}}
\put(90,130){\circle*{5}}


\end{picture}
\caption{}
\label{mainnandabesa}
\end{figure}

\bigskip

We can detemine all the possibilities for $\Gamma_{\Sigma}$ as follows: 
First, we classify lattice polygons in 
$$\Conv\left(
e_2,e_3,z_1
\right)\cup
\Conv\left(
e_2,e_3,z_2
\right)$$
which contain $e_2$, $e_3$ and $x_0+x_1$ (there exist exactly $23$ 
such polygons). Then, for any classified polygon $\Delta$, 
we detemine $\Gamma_{\Sigma}$ such that $F_{\sigma}=\Delta$ 
for a $3$-dimensional cone $\sigma\in\Sigma$. There exist exactly 
$13$ Gorenstein local toric Fano $3$-folds in this case. 
Table \ref{yokuwakaran} is the classified list for these 
$13$ Gorenstein local toric Fano $3$-folds. In this table, 
we give two data: (1) $\G(\Sigma)\setminus\{e_1,e_2,e_3\}$ and 
(2) the equation of $H_{\sigma}$ for any 
$3$-dimensional cone $\sigma\in\Sigma$. 
The following lemma is the key 
for this calculation. 

\begin{lem}
Let $X=X_{\Sigma}$ be a Gorenstein local toric Fano $3$-fold. 
For any $3$-dimensional cone $\sigma\in\Sigma$, 
$S_{\sigma}$ contains at least one of $e_1$, $e_2$ and $e_3$. 
\end{lem}

\begin{proof}
Let $aZ_1+bZ_2+cZ_3=1$ be the equation of $H_{\sigma}$ for a 
$3$-dimensional cone $\sigma\in\Sigma$. 
Obviously, $a\leq 1$, $b\leq 1$ and $c\leq 1$. 
One can easily confirm that 
at least one of $a$, $b$ and $c$ is equal to $1$.
\end{proof}

\begin{longtable}{|c||l|l|}
\caption{Gorenstein local toric Fano $3$-folds 
in the case (IIb)}\label{yokuwakaran} \\

\hline

\multicolumn{1}{|c||}{} & 
\multicolumn{1}{c|}{$\G(\Sigma)\setminus\{e_1,e_2,e_3\}$} & 
\multicolumn{1}{c|}{$H_{\sigma}$'s} \\
\hline
\hline

$(1)$ & 
$
\begin{pmatrix}
15\cr 10\cr 6
\end{pmatrix}
$.

& \begin{tabular}{@{}p{5.5cm}@{}}
$-Z_1+Z_2+Z_3=1$, 
$Z_1-2Z_2+Z_3=1$ and 
$Z_1+Z_2-4Z_3=1$.
\end{tabular}

\\
\hline

$(2)$ & 
$
\begin{pmatrix}
6\cr 4\cr 3
\end{pmatrix}
$.

& \begin{tabular}{@{}p{5.5cm}@{}}
$-Z_1+Z_2+Z_3=1$, 
$Z_1-2Z_2+Z_3=1$ and 
$Z_1+Z_2-3Z_3=1$.
\end{tabular}

\\
\hline
$(3)$ & 
$
\begin{pmatrix}
3\cr 2\cr 2
\end{pmatrix}
$.

& \begin{tabular}{@{}p{5.5cm}@{}}
$-Z_1+Z_2+Z_3=1$, 
$Z_1-2Z_2+Z_3=1$ and 
$Z_1+Z_2-2Z_3=1$.
\end{tabular}

\\
\hline

$(4)$ & 
$
\begin{pmatrix}
9\cr 6\cr 4
\end{pmatrix}
$ and 
$
\begin{pmatrix}
5\cr 4\cr 2
\end{pmatrix}
$.

& \begin{tabular}{@{}p{5.5cm}@{}}
$-Z_1+Z_2+Z_3=1$, 
$Z_1-2Z_2+Z_3=1$, 
$Z_1+Z_2-4Z_3=1$ and 
$Z_1-2Z_3=1$.
\end{tabular}

\\
\hline

$(5)$ & 
$
\begin{pmatrix}
5\cr 4\cr 2
\end{pmatrix}
$ and 
$
\begin{pmatrix}
3\cr 2\cr 2
\end{pmatrix}
$.

& \begin{tabular}{@{}p{5.5cm}@{}}
$-Z_1+Z_2+Z_3=1$, 
$Z_1-2Z_2+Z_3=1$, 
$Z_1+Z_2-4Z_3=1$ and 
$Z_1-Z_2=1$.
\end{tabular}

\\
\hline

$(6)$ & 
$
\begin{pmatrix}
4\cr 3\cr 2
\end{pmatrix}
$ and 
$
\begin{pmatrix}
3\cr 2\cr 2
\end{pmatrix}
$.

& \begin{tabular}{@{}p{5.5cm}@{}}
$-Z_1+Z_2+Z_3=1$, 
$Z_1-2Z_2+Z_3=1$, 
$Z_1+Z_2-3Z_3=1$ and 
$Z_1-Z_2=1$.
\end{tabular}

\\
\hline

$(7)$ & 
$
\begin{pmatrix}
3\cr 2\cr 2
\end{pmatrix}
$ and 
$
\begin{pmatrix}
2\cr 2\cr 1
\end{pmatrix}
$.

& \begin{tabular}{@{}p{5.5cm}@{}}
$-Z_1+Z_2+Z_3=1$, 
$Z_1-2Z_2+Z_3=1$, 
$Z_1+Z_2-3Z_3=1$ and 
$Z_1-Z_3=1$.
\end{tabular}

\\
\hline

$(8)$ & 
$
\begin{pmatrix}
5\cr 4\cr 2
\end{pmatrix}
$ and 
$
\begin{pmatrix}
2\cr 1\cr 2
\end{pmatrix}
$.

& \begin{tabular}{@{}p{5.5cm}@{}}
$-Z_1+Z_2+Z_3=1$, 
$Z_1-3Z_2+Z_3=1$, 
$Z_1+Z_2-4Z_3=1$ and 
$Z_1-Z_2=1$.
\end{tabular}

\\
\hline

$(9)$ & 
$
\begin{pmatrix}
4\cr 3\cr 2
\end{pmatrix}
$ and 
$
\begin{pmatrix}
2\cr 1\cr 2
\end{pmatrix}
$.

& \begin{tabular}{@{}p{5.5cm}@{}}
$-Z_1+Z_2+Z_3=1$, 
$Z_1-3Z_2+Z_3=1$, 
$Z_1+Z_2-3Z_3=1$ and 
$Z_1-Z_2=1$.
\end{tabular}

\\
\hline

$(10)$ & 
$
\begin{pmatrix}
6\cr 4\cr 3
\end{pmatrix}
$
, 
$
\begin{pmatrix}
5\cr 3\cr 2
\end{pmatrix}
$ 
and 
$
\begin{pmatrix}
3\cr 2\cr 1
\end{pmatrix}
$
.

& \begin{tabular}{@{}p{5.5cm}@{}}
$-Z_1+Z_2+Z_3=1$, 
$Z_1-2Z_2+Z_3=1$, 
$Z_1+Z_2-4Z_3=1$, 
$Z_1-2Z_3=1$ and 
$Z_2-Z_2=1$.
\end{tabular}
 \\
\hline

$(11)$ & 
$
\begin{pmatrix}
4\cr 3\cr 2
\end{pmatrix}
$
, 
$
\begin{pmatrix}
3\cr 2\cr 2
\end{pmatrix}
$ 
and 
$
\begin{pmatrix}
3\cr 2\cr 1
\end{pmatrix}
$
.

& \begin{tabular}{@{}p{5.5cm}@{}}
$-Z_1+Z_2+Z_3=1$, 
$Z_1-2Z_2+Z_3=1$, 
$Z_1+Z_2-4Z_3=1$, 
$Z_1-Z_2=1$ and 
$Z_2-Z_2=1$.
\end{tabular}
 \\
\hline

$(12)$ & 
$
\begin{pmatrix}
4\cr 3\cr 2
\end{pmatrix}
$
, 
$
\begin{pmatrix}
3\cr 2\cr 1
\end{pmatrix}
$ 
and 
$
\begin{pmatrix}
2\cr 1\cr 2
\end{pmatrix}
$
.

& \begin{tabular}{@{}p{5.5cm}@{}}
$-Z_1+Z_2+Z_3=1$, 
$Z_1-3Z_2+Z_3=1$, 
$Z_1+Z_2-4Z_3=1$, 
$Z_1-Z_2=1$ and 
$Z_2-Z_2=1$.
\end{tabular}
 \\
\hline

$(13)$ & 
$
\begin{pmatrix}
3\cr 2\cr 2
\end{pmatrix}
$
, 
$
\begin{pmatrix}
2\cr 1\cr 1
\end{pmatrix}
$ 
and 
$
\begin{pmatrix}
2\cr 1\cr 2
\end{pmatrix}
$
.

& \begin{tabular}{@{}p{5.5cm}@{}}
$-Z_1+Z_2+Z_3=1$, 
$Z_1-3Z_2+Z_3=1$, 
$Z_1+Z_2-3Z_3=1$, 
$Z_1-Z_3=1$ and 
$Z_1-Z_2=1$.
\end{tabular}
 \\
\hline

\end{longtable}

\bigskip

(III) Since 
$$x_1+x_2+e_3=
\begin{pmatrix}
3\cr 3\cr 3
\end{pmatrix}
=3x_0,$$
$\Conv(e_3,x_1,x_2)$ contains $x_0$ in its relative interior. 
So, $\Conv(e_3,x_1,x_2)$ is contained in the boundary of 
$\Gamma_{\Sigma}$. $\Conv(e_3,x_1,x_2)$ is contained in the 
hyperplane $(Z_3=1)\subset N_{\mathbb R}$. Therefore, 
by Proposition \ref{sakura}, 
$\Sigma$ is associated to a polytope $P\subset(Z_3=1)$ as 
in Example \ref{pkarasada}. 

\bigskip

(IV) In this case, 
$\Conv(x_1,x_2,x_3)$ is contained in $\Gamma_{\Sigma}$ by 
Proposition \ref{charfano}. So, 
$$x_1+x_2+x_3=
\begin{pmatrix}
4\cr 4\cr 4
\end{pmatrix}
=4x_0$$
implies that $x_0$ is contained in the interior of 
$\Gamma_{\Sigma}$. This is impossible. 
Therefore, the case (IV) does not occur. 

\bigskip

We end this section by summarizing the classification.

\renewcommand{\theenumi}{\roman{enumi}}

\begin{thm}\label{babyg}
Let $X=X_{\Sigma}$ be a Gorenstein local toric Fano $3$-fold and 
$X\neq{\mathbb A}^3$. 
Then, one of the following holds$:$
\begin{enumerate}
\item For some $m\in{\mathbb Z}_{>0}$,
$$\Gamma_{\Sigma}=\Conv\left(0,e_1,e_2,e_3,
\begin{pmatrix}
m\cr m\cr 1
\end{pmatrix}
\right)$$
$($see {\rm Proposition \ref{case1matome}}$)$.
\item $\Sigma$ is associated to a polytope $P\subset(Z_3=1)$ as 
in {\rm Example \ref{pkarasada}}. 
\item The other $13$ cases $($see {\rm Table \ref{yokuwakaran}}$)$.
\end{enumerate}
\end{thm}


\end{document}